# The Optimal Arbitrary-Proportional Finite-Set-Partitioning

Tiancheng Li

BISITE research group, University of Salamanca, Salamanca 37008, Spain; and
School of Mechanical Engineering, Northwestern Polytechnical University, Xi'an 710072, China.
tiancheng.li1985@gmail.com, t.c.li@mail.nwpu.edu.cn;
Webpage: https://sites.google.com/site/tianchengli85;

*Abstract*

This paper considers the arbitrary-proportional finite-set-partitioning problem which involves partitioning a finite set of size $N$ into $M$ subsets with respect to arbitrary nonnegative proportions $w^{(m)}, m = 1, 2, \ldots M$, where $N$ and $M$ are positive integers. This is the core art of many fundamental problems such as determining quotas for different individuals of different weights or sampling from a discrete-valued weighted sample set to get a new identically distributed but non-weighted sample set (e.g. the resampling needed in the particle filter). The challenge raises as the size of each subset, denoted as $N^{(m)}$, must be an integer while the unbiased expectation $Nw^{(m)}$ is often not, given that $\sum_{m=1}^{M} N^{(m)} = N, \sum_{m=1}^{M} w^{(m)} = 1$. To solve this problem, a metric (cost function) is defined on their discrepancies and correspondingly a solution is proposed to determine the sizes of each subsets, gaining the minimal cost. Theoretical proof and simulation demonstrations are provided to demonstrate the optimality of the scheme in the sense of the proposed metric.

*Keywords:* Finite-set-partitioning; discrete sampling; statistics analysis;

## I. Introduction & Problem Formulation

In the real world, we often come across a type of allocation/collection problems that need to allocate, or conversely to collect, $N$ equivalent goods to/from $M$ individuals with respect to the nonnegative proportion $w^{(m)}, m = 1,2,\ldots M$, where $N$ and $M$ are positive integers and the sum of the proportions is one. The unbiased expectation of the number of goods allocated to/collected from individual $m$, denoted as $N^{(m)}$ for which $\sum_{m=1}^{M} N^{(m)} = N$, is

$$E(N^{(m)}|w^{(m)},N) = Nw^{(m)} \tag{1}$$

This is the core art of many fundamental problems, which is referred to as arbitrary-proportional finite-set-partitioning in this paper, involving partitioning the finite set of size $N$ into $M$ subsets with respect to the nonnegative proportions $w^{(m)}, m = 1,2,\ldots M$, where $\{N, M\} \in \mathbb{N}^+, \sum_{m=1}^{M} w^{(m)} = 1$. The difficulty of this problem lies in the situation in which the size of the subset $N^{(m)}$ must be an integer while the unbiased expectation $Nw^{(m)}$ is often not otherwise we can straightforwardly have $N^{(m)} = Nw^{(m)}$. Then, how to determine the nonnegative integer sequence $\{N^{(m)}\}_{m=1}^{M}$ *optimally* so that the bias (to be precisely defined) is minimal while satisfying all the practical constraints?

To solve this problem, a metric is required firstly to define the bias which can be then used as the cost function for optimization. It is interesting/reasonable to care about the discrepancies between the sizes of the subsets and the unbiased expectation, namely $N^{(m)} - Nw^{(m)}$ for $m = 1,2,\ldots M$. The larger the discrepancies are, the worse (more biased) the solution is. So far, the problem can be completely modelled as

$$\{N^{(m)}\}_{m=1}^{M} = \operatorname{argmin}_{\{N^{(m)}\}_{m=1}^{M}} \frac{1}{M} \sum_{m=1}^{M} \|N^{(m)} - Nw^{(m)}\| \tag{2}$$

subject to $\begin{cases} \{N, M\} \in \mathbb{N}^+ \\ \{w^{(m)}\}_{m=1}^{M} \in \mathbb{R}_0^+; \sum_{m=1}^{M} w^{(m)} = 1 \\ \{N^{(m)}\}_{m=1}^{M} \in \mathbb{N}_0^+; \sum_{m=1}^{M} N^{(m)} = N \end{cases}$

where $\|\cdot\|$ is a distance/metric to be specified, $\mathbb{N}_0^+$ is the set of nonnegative integers, $\mathbb{N}^+$ is the set of positive integers and $\mathbb{R}_0^+$ is the set of nonnegative real numbers.

## II. The metric and the solution

In this paper, we define the metric required in (2) on the second moment of the discrepancy, referred to as the mean square error (MSE) which is similar to the statistical variance, as

$$MSE = \frac{1}{M} \sum_{m=1}^{M} (N^{(m)} - Nw^{(m)})^2 \tag{3}$$

Obviously, the MSE provides a consistent and efficient metric to quantify/measure the overall discrepancy between the unbiased expectation and the result obtained by a finite-set-partitioning solution. To obtain the optimal result of the minimal bias, it is desirable to determine the nonnegative integral size of each subset so that MSE is the least/minimum. In the following, we will present such a deterministic finite-set-partitioning scheme, referred to the least MSE (LMSE) finite-set-partitioning.

The LMSE finite-set-partitioning procedure consists of two parts. First, determine the initial subset sizes $N^{(m)} = \text{Floor}(N \times w^{(m)})$ for $m = 1, 2, \ldots, M$ and we obtain $L = \sum_{n=1}^{M} N^{(m)}$, where $\text{Floor}(x)$ gives the largest integer not exceeding $x$. Secondly, rank the proportion residuals $w^{(m)} - N^{(m)}/N$ for $m = 1, 2, \ldots, M$ and find the largest $N - L$ residuals for each of which add the corresponding subset size by one. Finally, we have $\sum_{n=1}^{M} N^{(m)} = L + N - L = N$. The scheme can be described as given in Algorithm 1, where the function $\text{TopRank}_s[S]$ gives a subset that contains the largest $s$ elements in the set of $S$.

Obviously, the proposed procedure as shown in Algorithm 1 is based on the metric of mean absolute error (MAE) $MAE = \frac{1}{M}\sum_{m=1}^{M}|N^{(m)} - Nw^{(m)}|$ to determine the sizes of the subsets. What follows will further proof the LMSE attributes of the algorithm.

**Algorithm 1** LMSE finite-set-partitioning

---

Input: $\left[\{w^{(m)}\}_{m=1}^{M}, N\right]$

Output: $\left[\{N^{(m)}\}_{m=1}^{M}\right]$

Procedure:

$L = 0$

FOR $m = 1: M$

    $N^{(m)} = \text{Floor}(N \times w^{(m)})$

    $\hat{w}^{(m)} = w^{(m)} - N^{(m)}/N$

    $L = L + N^{(m)}$

END

$\{\hat{w}^{(l)}\}_{l=1}^{N-L} = \text{TopRank}_{N-L}\left[\{\hat{w}^{(m)}\}_{m=1}^{M}\right]$

FOR $m = 1: M$

    IF $\hat{w}^{(m)} \in \{\hat{w}^{(l)}\}_{l=1}^{N-L}$

        $N^{(m)} = N^{(m)} + 1$

    END

END

---

**Theory 1**. The output of Algorithm 1 satisfies that $\forall\ 1 < m < M$, $|N^{(m)} - Nw^{(m)}| < 1$.

**Proof.** As shown in Algorithm 1, it is straightforward to know that if $\hat{w}^{(m)} \in \{\hat{w}^{(l)}\}_{l=1}^{N-n}$, $N^{(m)} = \text{Floor}(N \times w^{(m)}) + 1$ otherwise $N^{(m)} = \text{Floor}(N \times w^{(m)})$. For both cases, we have $\forall\ 1 < m < M$, $|N^{(m)} - Nw^{(m)}| < 1$ as stated.

**Theory 2**. A necessary condition for the LMSE of Eq. (3) is that $\forall\ 1 < m < M$, $|N^{(m)} - Nw^{(m)}| < 1$.

**Proof.** Assuming a subset of size $N^{(q)}$ satisfies $(N^{(q)} - Nw^{(q)}) \geq 1$, there must exist another subset $p: (N^{(p)} - Nw^{(p)}) < 0$ due to the overall constraint $\sum_{m=1}^{M}(N^{(m)} - Nw^{(m)}) = 0$.

If we change the result to be $N_{new}^{(p)} = N^{(p)} + 1$ and $N_{new}^{(q)} = N^{(q)} - 1$, while $N_{new}^{(i)} = N^{(i)}, i \neq p, q$, then we have the MSE change as follows

$$MSE_{new} - MSE = \frac{1}{M} \sum_{i=\{p,q\}} \left(N_{new}^{(i)} - Nw^{(i)}\right)^2 - \frac{1}{M} \sum_{i=\{p,q\}} \left(N^{(i)} - Nw^{(i)}\right)^2$$

$$= \frac{1}{M} \sum_{i=\{p,q\}} \left(\left(N_{new}^{(i)}\right)^2 - \left(N^{(i)}\right)^2 + 2N^{(i)}Nw^{(i)} - 2N_{new}^{(i)}Nw^{(i)}\right)$$

$$= \frac{1}{M}\left(\left(N^{(p)} + 1\right)^2 - \left(N^{(p)}\right)^2 + 2N^{(p)}Nw^{(p)} - 2\left(N^{(p)} + 1\right)Nw^{(p)}\right)$$

$$+ \frac{1}{M}\left(\left(N^{(q)} - 1\right)^2 - \left(N^{(q)}\right)^2 + 2N^{(q)}Nw^{(q)} - 2\left(N^{(q)} - 1\right)Nw^{(q)}\right)$$

$$= \frac{2}{M}\left(1 + \left(N^{(p)} - Nw^{(p)}\right) - \left(N^{(q)} - Nw^{(q)}\right)\right)$$

$$< \frac{2}{M}(1 + 0 - 1) = 0 \tag{4}$$

This indicates that $\left(N^{(q)} - Nw^{(q)}\right) \geq 1$ cannot be a part of the LMSE result. The similar proof and results hold to the case of $\left(N^{(q)} - Nw^{(q)}\right) \leq -1$ and the detailed process of the proof is omitted here. Therefore, we conclude that the LMSE of (3) requires $\forall\ 1 < m < M,\ \left|N^{(m)} - Nw^{(m)}\right| < 1$.

**Theory 3**. Algorithm 1 achieves the LMSE solution.

Proof. Theories 1 and 2 indicate that the output of Algorithm 1 and the LMSE result are possibly the same or they can become each other through adjusting of $N^{(m)}$ between two particles one/several time(s). If we can proof that any change on $N^{(m)}$ given by Algorithm 1 will cause an increase of SE, we will be certain that the output of Algorithm 1 is the LMSE. Without loss of generality, assume the change occurs on the subsets $p$ and $q$ of the size $N^{(p)}$ and $N^{(q)}$ respectively as determined by Algorithm 1. For $N_{new}^{(p)} = N^{(p)} + l$ and $N_{new}^{(q)} = N^{(q)} - l$, where $1 \leq l \leq \min(N - L, N^{(q)})$, while $N_{new}^{(i)} = N^{(i)}, i \neq p, q$, we have

$$MSE_{new} - MSE = \frac{1}{M}\sum_{i=\{p,q\}}\left(N_{new}^{(i)} - Nw^{(i)}\right)^2 - \sum_{i=\{p,q\}}\left(N^{(i)} - Nw^{(i)}\right)^2$$

$$= \frac{1}{M}\sum_{i=\{p,q\}}\left(\left(N_{new}^{(i)}\right)^2 - \left(N^{(i)}\right)^2 + 2N^{(i)}Nw^{(i)} - 2N_{new}^{(i)}Nw^{(i)}\right)$$

$$= \frac{1}{M}\binom{\left(N^{(p)} + l\right)^2 - \left(N^{(p)}\right)^2 + 2N^{(p)}Nw^{(p)} - 2\left(N^{(p)} + l\right)Nw^{(p)}}{+\left(N^{(q)} - l\right)^2 - \left(N^{(q)}\right)^2 + 2N^{(q)}Nw^{(q)} - 2\left(N^{(q)} - l\right)Nw^{(q)}}$$

$$= \frac{2l}{M}\left(l + \left(N^{(p)} - Nw^{(p)}\right) - \left(N^{(q)} - Nw^{(q)}\right)\right) \tag{5}$$

The proposed method given in Algorithm 1 satisfies the following bounds

$$N^{(i)} - Nw^{(i)} = \begin{cases} -\widehat{w}^{(i)} & N^{(i)} \leq Nw^{(i)} \\ 1 - \widehat{w}^{(i)} & N^{(i)} > Nw^{(i)} \end{cases} \tag{6}$$

and $0 < \widehat{w}^{(p)}, \widehat{w}^{(q)} < 1 \leq l$. Specifically, if $N^{(p)} \leq Nw^{(p)}$ and $N^{(q)} > Nw^{(q)}$, we have $\widehat{w}^{(q)} > \widehat{w}^{(p)}$. Then, Eq. (5) will go to the following four cases:

1) If $N^{(p)} \leq Nw^{(p)}, N^{(q)} \leq Nw^{(q)}$, (5) will reduce to

$$MSE_{new} - MSE = \frac{2l}{M}\left(l - \widehat{w}^{(p)} + \widehat{w}^{(q)}\right) > 0$$

2) If $N^{(p)} \leq Nw^{(p)}, N^{(q)} > Nw^{(q)}$, we have $\widehat{w}^{(q)} > \widehat{w}^{(p)}$ and (5) will reduce to

$$MSE_{new} - MSE = \frac{2l}{M}\left(l - 1 + \widehat{w}^{(q)} - \widehat{w}^{(p)}\right) > 0$$

3) If $N^{(p)} > Nw^{(p)}, N^{(q)} \leq Nw^{(q)}$, (5) will reduce to

$$MSE_{new} - MSE = \frac{2l}{M}\left(l - \widehat{w}^{(p)} + \widehat{w}^{(q)}\right) > 0$$

4) If $N^{(p)} > Nw^{(p)}, N^{(q)} > Nw^{(q)}$, (5) will reduce to

$$MSE_{new} - MSE = \frac{2l}{M}\left(l - \widehat{w}^{(p)} + \widehat{w}^{(q)}\right) > 0$$

It is shown that all possible changes on the output of Algorithm 1 will lead to an increase of the MSE. This demonstrates that Algorithm 1 achieves the least MSE.

## III. LMSE sampling

The proposed solution can serve as a resampling scheme for the particle filter which is to sample a new set of equally weighted particles from a set of particles that are of different weights, in order to reduce the weight degeneracy problem [2]. Here, **M** is the number of particles before resampling, **N** is the number of new equal-weighted particles obtained after resampling and $N^{(m)}$ is the number of times that particle **m** is resampled. Here, the MSE metric can be called as *sampling variance* (SV) and the proposed resampling method is called minimum-sampling-variance (MSV) resampling. To note, a smaller SV means a better identical distribution attribute of the resampled particles as compared to the original particles before resampling. Identical distribution is a critical principle for the discrete sampling.

To demonstrate the LMSE/MSV performance of the proposed solution, the MSV resampling method is employed in the general sampling-importance-resampling (SIR) filter on a one-dimension filtering model with the state transition equation and the observation equation respectively given as follows

$$x_t = 1 + \sin(w\pi t) + \phi_1 x_{t-1} + u_t \tag{7}$$

$$y_t = \begin{cases} \phi_2 x_t^2 + v_t & t \leq 30 \\ \phi_3 x_t - 2 + v_t & t > 30 \end{cases} \tag{8}$$

where $x_t, y_t$ are respective the state and measurement at iteration step $t$, the scale parameters $\omega = 4e - 2$, $\phi_1 = 0.5, \phi_2 = 0.2$ and $\phi_3 = 0.5$, the process noise $u_t$ is a Gamma $\mathcal{G}a(3,2)$ random variable and the observation noise is Gaussian $v_t \sim \mathcal{N}(0,1)$. To note here, the observation model (8) uses a relative large observation noise variance in which the general SIR can work well without using the unscented Kalman filter or even others as the sampling proposal [1].

Traditional resampling methods include the multinomial resampling, residual resampling, systematic resampling and residual systematic resampling (RSR) (their pseudo-codes are given in [2] and the MATLAB codes are available online: https://sites.google.com/site/tianchengli85/matlab-codes/resampling-methods) are also realized for comparison. All these resampling methods share the same attribute that they output the same number of particles and all the resampled particles are equally weighted. In particular, the multinomial resampling generates (random) independent identical distributed particles while the other three resampling schemes employ different degrees of deterministic sampling operations. Comparably, the MSV resampling uses no random number and is purely deterministic sampling.

Since we are only interested in the SV performance, different resampling methods are employed on exactly the same particles population that are generated by the same SIR filter at each time-step. The filter uses the same, a single, resampling for filtering and 100 particles are used. The state and estimates against time are given in the upper subfigure of Fig. 1 while the SV obtained by different resampling methods are given in the middle and the bottom (zoomed in) sub-figures. The result explicitly show that the MSV resampling approach achieves the smallest SV among all. Reasonably, the performance of the systematic resampling and the RSR resampling are very similar, the SV obtained by the multinomial resampling is the highest while the SV obtained by the residual resampling is in the middle of them.

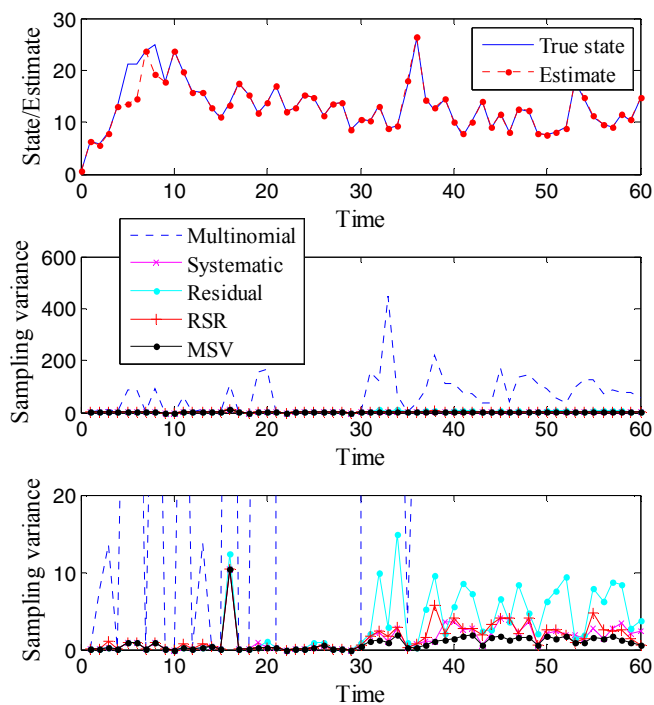

**Fig. 1** State estimation and the sampling variance against time

## IV. Conclusion

This paper concerns the arbitrary-proportional discrete-valued finite-set-partitioning problem, namely partitioning a finite set of size $N$ into $M$ subsets in proportion as specified. Under the constraint that each subset must be of an integer size, the optimal solution is established that is of the least mean square error as compared to the unbiased expectation determined by the proportions. The proposed scheme can be used for resampling in the particle filter, namely the minimum-sampling-variance resampling which is a deterministic sampling scheme and will always generate the particle set of the minimum-sampling-variance attribute. The optimality of the proposed solution is demonstrated in theory and in real-data via simulations.